\documentclass[11pt]{amsart}

\makeatletter
\usepackage{amssymb}
\usepackage{latexsym}
\usepackage{amsbsy}
\usepackage{amsfonts}

\def\marginpar#1{\ignorespaces}

\newtheorem{theorem}[equation]{Theorem}
\newtheorem{proposition}[equation]{Proposition}
\newtheorem{lemma}[equation]{Lemma}
\newtheorem{corollary}[equation]{Corollary}

\theoremstyle{definition}
\newtheorem{remark}[equation]{Remark}

\newtheorem{example}[equation]{Example}

\numberwithin{equation}{section}

\def\AArm{\fam0 \rm}%
\newdimen\AAdi%
\newbox\AAbo%
\def\AAk#1#2{\setbox\AAbo=\hbox{#2}\AAdi=\wd\AAbo\kern#1\AAdi{}}%

\newcommand{\BBone}{{\ensuremath{{\AArm 1\AAk{-.8}{I}I}}}}

\def\eqref#1{(\ref{#1})}
\def\eqlabel#1{\def\@currentlabel{#1}}

\def\formula#1{\def\@tempa{#1}\let\@tempb\theequation\def\theequation{%
\hbox{#1}}\def\@currentlabel{(\theequation)}$$}
\def\endformula{\leqno\hbox{(\@tempa)}$$\@ignoretrue\let\theequation\@tempb}

\def\given{\hskip5\p@\relax\vrule\@width.4\p@\hskip5\p@\relax}

\newcommand{\open}[1]{%
\par\normalfont\topsep6\p@\@plus6\p@\trivlist\item[\hskip\labelsep\itshape#1%
\@addpunct{.}]\ignorespaces}

\DeclareRobustCommand{\close}[1]{%
  \ifmmode 
  \else \leavevmode\unskip\penalty9999 \hbox{}\nobreak\hfill
  \fi
  \quad\hbox{$#1$}}

\newlength{\toskip}

\def\DD{\mathcal D}

\def\EE{\mathcal E}

\def\LL{\mathcal L}

\def\<{\langle}
\def\>{\rangle}

\def \R {{\mathbb R}}

\def \N {{\mathbb N}}
\def \D {{\mathbb D}}

\def \Var {\textrm{Var}}
\def \Ent {\textrm{Ent}}

\makeatother

\begin{document}
\date{\today}

\title[Talagrand's transportation inequality]{A note on Talagrand's transportation inequality and logarithmic Sobolev inequality}

 \author[P. Cattiaux]{\textbf{\quad {Patrick} Cattiaux $^{\spadesuit}$ \, \, }}
\address{{\bf {Patrick} CATTIAUX},\\ Institut de Math\'ematiques de Toulouse. CNRS UMR 5219. \\
Universit\'e Paul Sabatier, Laboratoire
de Statistique et Probabilit\'es,
\\ 118 route
de Narbonne, F-31062 Toulouse cedex 09.} \email{cattiaux@math.univ-toulouse.fr}

 \author[A. Guillin]{\textbf{\quad {Arnaud} Guillin $^{\diamondsuit}$}}
\address{{\bf {Arnaud} GUILLIN},\\ Laboratoire de Math\'ematiques, CNRS UMR 6620, Universit\'e Blaise Pascal, avenue des Landais 63177 Aubi\`ere.}
\email{guillin@math.univ-bpclermont.fr}

 \author[L-M. Wu]{\textbf{\quad {Li-Ming} Wu $^{\diamondsuit}$}}
\address{{\bf {Li-Ming} WU},\\ Laboratoire de Math\'ematiques, CNRS UMR 6620, Universit\'e Blaise Pascal, avenue des Landais 63177 Aubi\`ere.}
\email{wuliming@math.univ-bpclermont.fr}

\maketitle
 \begin{center}

 \textsc{$^{\spadesuit}$  Universit\'e de Toulouse}
\smallskip

\textsc{$^{\diamondsuit}$ Universit\'e Blaise Pascal}
\smallskip

 \end{center}

\begin{abstract}
We give by simple arguments sufficient conditions, so called Lyapunov conditions, for Talagrand's transportation information inequality and for the logarithmic Sobolev inequality. Those sufficient conditions work even in the case where the Bakry-Emery curvature is not lower bounded. Several new examples are provided.
\end{abstract}
\bigskip

\bigskip

\bigskip

\section{\bf Introduction and main results.}\label{Intro}

Transportation cost information inequalities have been recently deeply studied,  especially for
their connection with the concentration of measure phenomenon, or for deviation inequalities for
Markov processes (see \cite{Led01,GLWY07}). In particular, Talagrand \cite{Tal96a} establishes the
so-called $T_2$ inequality (or Talagrand's transportation inequality, or $W_2H$ inequality) for
the Gaussian measure, establishing thus Gaussian dimension free concentration of measure. But
before going further in the numerous results around these inequalities, let us present the object
under study.

Given a metric space $(E,d)$ equipped with its Borel $\sigma$ field, and $1\leq p < +\infty$, the
$L^p$ Wasserstein distance between two probability measures $\mu$ and $\nu$ on $E$ is defined as
\begin{equation}\label{wass}
W_p(\mu,\nu):=\left(\inf_{\pi}\int_{E\times E} d^p(x,y)~\pi(dx,dy)\right)^{1/p}
\end{equation}
where the infimum runs over all coupling $\pi$ of $(\mu,\nu)$, see Villani \cite{Vill03} for an
extensive study of such quantities.

A probability measure $\mu$ is then said to satisfy the transportation-entropy inequality
$W_pH(C)$,  where $C>0$ is some constant, if for all probability measure $\nu$
\begin{equation}\label{Tp}
W_p(\nu,\mu)\le \sqrt{2C~H(\nu|\mu)}
\end{equation}
where $H(\nu|\mu$) is the Kullback-Leibler information, or relative entropy, of $\nu$ with respect to $\mu$:
\begin{equation}
H(\nu|\mu):=\left\{\begin{array}{ll}\int\log\left(\frac{d\nu}{d\mu}\right)d\nu&{\rm if}~\nu\ll\mu\\
+\infty&{\rm otherwise.}\end{array}\right.
\end{equation}
Marton \cite{Mar96} has first shown how $W_1H$ inequality implies Gaussian concentration  of
measure and Talagrand, via a tensorization argument, established that the standard Gaussian
measure, in any dimension, satisfies $W_2H(C)$ with the sharp constant $C=1$.

However, if $W_1H$ is completely characterized via a practical
Gaussian integrability  criterion (see \cite{DGW03,BV05}),  $W_2H$
is much more difficult to describe. Nevertheless several equivalent
beautiful conditions are known.

\begin{theorem}
The following conditions are equivalent
\begin{enumerate}
\item $\mu$ satisfies $W_2H(C)$ for some constant $C>0$.
\item For any bounded and measurable function $f$ with $\mu(f)=0$, defining the inf-convolution
$$Qf(x)=\inf_{y\in E}\{f(y)+d^2(x,y)\},$$ we have
\begin{equation}
\label{BG} \int e^{\frac 1{2C}Qf}d\mu\le 1.
\end{equation}
\item There exist $a$, $r_0$, $b$ such that for all $n$ all measurable $A\subset E^n$, with
$\mu^{\otimes n}(A)\ge 1/2$, the probability measure $\mu^{\otimes n}$ satisfies
\begin{equation}
\label{gozadim} \mu^{\otimes n}(A^r)\ge 1-b~e^{-a(r-r_0)^2}
\end{equation}
where $A^r=\{x\in E^n;~\exists y\in A, \sum_1^nd^2(x_i,y_i)\le r^2\}.$
\end{enumerate}
\end{theorem}
$(1) \Leftrightarrow (2)$ was proved in the seminal paper by
Bobkov-G\"otze \cite{BG99}, and $(1) \Leftrightarrow (3)$ very
recently by Gozlan \cite{Goz08}. Hence we have the beautiful
characterization, $W_2H$ is nothing else than a dimension free
Gaussian concentration for the product measure. Note also that
Gozlan-L\'eonard \cite{GL07} established another criterion as a
large deviation upper bound. One point is however important to
remark: if these various characterizations have nice implications
(concentration, deviation,...), it is rather difficult to directly
use them to prove a $W_2H$ inequality.

The first step towards practical criterion was done by Otto-Villani \cite{OVill00}, soon followed
by Bobkov-Gentil-Ledoux \cite{BGL01}, who established that if $\mu$ satisfies a logarithmic
Sobolev inequality, then $\mu$ satisfies $W_2H$ (note that many explicit sufficient conditions for
log-Sobolev inequalities are now known). To be more precise, let us present our framework.

Throughout this paper $E$ is a complete and connected Riemannian manifold of finite dimension, $d$
the geodesic distance, and $dx$ the volume measure. $\mu(dx)=e^{-V(x)}dx/Z$ is the Boltzmann
measure  with $V\in C^2$ and $Z=\int e^{-V}dx<+\infty$. If the logarithmic Sobolev inequality
$LSI(C)$ is verified, i.e. for all locally lispchitz $g$
\begin{equation}\label{LSI}
\Ent_\mu(g^2):=\int g^2\log\left(\frac{g^2}{\int g^2d\mu}\right)d\mu\le 2C \int |\nabla g|^2d\mu
\end{equation}
then $\mu$ satisfies also $W_2H(C)$. The proof of Otto-Villani
\cite{OVill00} relies on a dynamical approach,  namely to derive the
Wasserstein distance between $\nu_t$ and $\nu_{t+s}$ when $\nu_t$ is
the dynamical transport leading from $\nu$ to $\mu$, whereas
Bobkov-Gentil-Ledoux \cite{BGL01} apply the hypercontractivity of
the Hamilton semigroup, leading to an Herbst's like argument to
derive $W_2H$.

It is only a few years ago that the two first authors \cite{CG05}
succeeded in  proving that $W_2H$ is strictly weaker than $LSI$,
providing an example in one dimension of a measure (with unbounded
curvature) satisfying $W_2H$ but not $LSI$. Their method is a
refinement of the argument of Bobkov-Gentil-Ledoux \cite{BGL01}:
indeed, a full $LSI$ is too strong to give $W_2H$, a $LSI$ for a
restricted class of functions is sufficient. They were however only
able to give an explicit sufficient condition in dimension one for
this restricted inequality. We will give here a Lyapunov condition
ensuring that this restricted logarithmic Sobolev inequality holds,
and thus $W_2H$ too. We will also show that if the Bakry-Emery
curvature $Ric + Hess_V$  is lower bounded then the same condition
implies $LSI$.

Consider the $\mu$-symmetric operator $\LL=\Delta -\nabla V.\nabla$ on E.  A Lyapunov condition is
of the form: there exists $W\ge 1$ and $r,b>0$ such that for some positive function $\phi$
\begin{equation}\label{genlyap}
\LL W\le -\phi W+b\BBone_{B(x_0,r)}.
\end{equation}
Such Lyapunov conditions have been used a lot both in discrete and continuous time case to  study
the speed of convergence towards the invariant measure of the associated semigroup under various
norms, see \cite{MT, DMT,DFG}. The deep connection between such conditions and various form of
functional inequalities have been recently studied by the authors (and coauthors). For example, if
$\phi$ is constant, it is shown in \cite{BBCG} that the Lyapunov condition implies both a
Poincar\'e inequality and a Cheeger inequality (with some slight additional assumptions on $W$).
If $\phi:=\phi(W)$ and $\phi$ is sub-linear then optimal weak Poincar\'e or isoperimetric
inequalities can be established, see \cite{BCG,CGGR}. Finally if $\phi:=\phi(W)$ is super-linear,
then it is shown to imply super Poincar\'e inequalities \cite{CGWW}, and thus various $F$-Sobolev
inequalities including logarithmic Sobolev inequalities.

Their implications in transportation cost inequalities were up to now not explored. It is the
purpose of this short note.

Here is our main result:

\begin{theorem}\label{thm1}
Let $\mu$ be a Boltzmann measure.\\
1)  Suppose  that there exists a $C^2$-function $W:E\to[1,\infty[$,
some point $x_0$ and constants  $b, c>0$ such that
\begin{equation}\label{lyapT2}
\LL W\le (-c d^2(x,x_0)+b)~W, \ x\in E
\end{equation}
or more generally there exists some nonnegative locally Lipschitzian
function $U$ ($=\log W$) such that in the distribution sense (see
the remark below),
\begin{equation}\label{lyapT2b}
\LL U +|\nabla U|^2\le -c d^2(x,x_0)+b
\end{equation}
then $W_2H(C)$ holds for some constant $C>0$.

 2) Under the Lyapunov condition (\ref{lyapT2}), suppose moreover that  $Hess(V)+Ric\ge K Id$ for
 some $K\le 0$ (in the sense of matrix). Then the logarithmic Sobolev inequality \eqref{LSI} holds.
\end{theorem}

\begin{remark}\label{rem1}
\begin{enumerate}
\item In both cases, it is of course possible to track all the constants  involved to get an upper
bound of the constant of $W_2H(C)$ inequality and of the logarithmic Sobolev inequality, as will
be seen from the proof. One will also remark that contrary to \cite{BCG,CGGR,CGWW}, we will not
use localization technique, constants are thus easier to derive.

\item If $U=\log W\in C^2$, then $\LL U +|\nabla U|^2=-\LL W/W$ so that
(\ref{lyapT2}) and (\ref{lyapT2b}) are equivalent. The condition
(\ref{lyapT2b}) in the distribution sense means that for any $h\in
C_0^\infty(E)$ (the space of infinitely differentiable functions
with compact support) such that $h\ge 0$,
$$\aligned
\int (\LL U +|\nabla U|^2) h dx&:=\int U \Delta h dx +\int\left(-
\nabla V\cdot \nabla U + |\nabla U|^2\right) h dx\\
& \le \int (-c d^2(x,x_0)+b) h d\mu. \endaligned$$

\item The Lyapunov condition (\ref{lyapT2}) implies that there exists $r_0>0$ and $b',\lambda>0$, such that
$$\LL W\le -\lambda W+b'\BBone_{B(x_0,r_0)}$$
so that, by \cite{BBCG}, $\mu$ satisfies a Poincar\'e inequality.

\end{enumerate}
\end{remark}

This paper is organized as follows. In the next section we present
several corollaries and examples for showing the usefulness and
sharpness of the Lyapunov condition (\ref{lyapT2b}). The very simple
proof of Theorem \ref{thm1} is given in Section 3. And in the last
section we combine the above-tangent lemma and the Lyapunov function
method to yield the LSI in the unbounded curvature case.


\section{Corollaries and examples}

\subsection*{Some practical conditions}

From Theorem \ref{thm1}, one easily deduces

\begin{corollary}\label{cor1}
Suppose that $\mu$ is a Boltzmann measure on $E=\R^d$. Let $x\cdot y$
and $|x|=\sqrt{x\cdot x}$ be the Euclidean inner product and norm, respectively.\\
1) If one of the following conditions
\begin{equation}\label{kustr}
\exists a<1, R, c>0, {\rm such~that~if~ }\ |x|>R,\qquad (1-a)|\nabla
V|^2-\Delta V\ge c~|x|^2
\end{equation}
or
\begin{equation}\label{simpl}
\exists  R,c>0, {\rm such~that~}\forall |x|>R, \qquad x\cdot\nabla
V(x)\ge c~|x|^2
\end{equation}
is satisfied, then $W_2H$ holds.\\
2) Under the same conditions, suppose moreover that  $Hess(V)\ge K
Id$ then a logarithmic Sobolev inequality (LSI in short) holds.
\end{corollary}

\begin{proof} Under (\ref{kustr}), one takes $W=e^{aV}$; and under
(\ref{simpl}) one choose $W=e^{a|x|^2}$ with $0<a<c/2$. One sees
that condition (\ref{lyapT2}) is satisfied in both case.
\end{proof}

\begin{remark}
\begin{enumerate}
\item Condition (\ref{kustr}) is of course reminiscent to the Kusuoka-Stroock condition for
logarithmic Sobolev inequality (replace $d^2$ by $V$). On the real
line, it implies the condition of \cite[Prop. 5.5]{CG05}. \item
Gozlan \cite[Prop. 3.9 and Theorem 4.8]{Goz08} proves $W_2H$ on
$\R^d$ under the condition
$$\liminf_{|x|\to\infty} \sum_{i=1}^d \left[\frac14 \left(\frac{\partial V}{\partial x_i}\right)^2-
\frac{\partial^2 V}{\partial x_i^2}\right]\frac{1}{1+x_i^2}\ge c$$
for some positive $c$, using weighted Poincar\'e inequality. Note that this condition is in general not
comparable to ours,  for the terms in the sum can be negative, and
also for we have more freedom with the choice of $a$ (limited to
$3/4$ in Gozlan's method). Whether this condition can be retrieved from a right choice of $W$ in (\ref{lyapT2}) seems unlikely. We will however simply show how to retrieve (and generalize) Gozlan's like conditions in the last section.

\item Condition (\ref{simpl}) may also be
compared with condition (1.7) in \cite{BBCG}: $x\cdot\nabla V(x)\ge
c~d(x,x_0)$ which implies Poincar\'e inequality.
\end{enumerate}
\end{remark}

\subsection*{Comparison with Wang's criterion}
 Wang's criterion for  LSI says the following: if $Hess_V + Ric \ge K
 Id$ with $K\le 0$ and $$\int e^{(|K|/2 + \varepsilon) d^2(x,x_0)}
 d\mu(x)<+\infty,$$ then $\mu=e^{-V} dx/C$ satisfies the LSI. We give
 now an example for which the previous criterion does not apply, but
 ours does.

\begin{example} {\rm Let $E=\R^2$ and $V(x,y)= r^2 g(\theta)$ for
all $r:=\sqrt{x^2+y^2}\ge 1$ (and $V\in C^\infty(\R^2)$), where
$(r,\theta)$ is the polar coordinates system and $g(\theta)=2+\sin
(k \theta)$ $(k\in \N^*)$ for all $\theta\in S^1\equiv [0,2\pi]$. We
have for $r>1$,
$$
(x,y)\cdot \nabla V(x,y)=r\partial_r V= 2r^2 g(\theta)\ge 2 r^2
$$
i.e., the condition (\ref{simpl}) is satisfied. Moreover $Hess_V$ is
bounded. Thus by Corollary \ref{cor1}, $\mu=e^{-V}dxdy/C$ satisfies
the LSI.

However Wang's integrability condition is not satisfied for large $k$. Indeed $\Delta V = 4
g(\theta) + g^{\prime\prime}(\theta)=8+(4-k^2)\sin\theta$, then the smallest eigenvalue
$\lambda_{\min}$ of $Hess_V$ satisfies
$$
\lambda_{\min}\le \frac 12 tr(Hess_V)=\frac 12\Delta V =4+ (2-k^2/2)
\sin(k\theta).
$$
Then the largest constant $K$ so that $Hess_V\ge K Id$ in the case
$k\ge 2$ satisfies
$$
K\le 6-k^2/2.
$$
When $k\ge 4$, $K/2\le 3-k^2/4\le -1$ and Wang's integrability
condition is not satisfied for $\int e^{r^2} d\mu=+\infty$. In other
words Wang's criteria does not apply for this example once $k\ge
4$.}
\end{example}

\subsection*{Riemannian manifold with unbounded curvature} Let $E$ be a
$d-$dimensional $(d\ge 2)$ connected complete Riemannian manifold
with
\begin{equation}\label{cor2a}
Ric_x \ge - (c + \sigma^2 d^2(x,x_0)), \ x\in E
\end{equation}
for some constants $c,\sigma>0$, where $x_0$ is some fixed point
$x_0$.  Let $V\in C^2(E)$ such that
\begin{equation}\label{cor2b}
\<\nabla d(x,x_0), \nabla V\> \ge \delta d(x,x_0)-k\ \text{ outside
of $cut(x_0)$ for some constants}\ \delta,k>0.
\end{equation}
Here $cut(x_0)$ denotes the the cut-locus of $x_0$.
\begin{corollary}\label{cor2} Assume  (\ref{cor2a}) and (\ref{cor2b}).
If $\delta>\sigma \sqrt{d-1}$, then $\mu=e^{-V} dx/C$ satisfies
$W_2H(C)$.
\end{corollary}

\begin{remark} Assume that $Hess_V\ge \delta$. Pick some $x\notin
cut(x_0)$, and denote by $U$ the unit tangent vector along the minimal geodesic $(x_s)_{0\le s\le
d(x,x_0)}$ from $x_0$ to $x$, we have
$$
\<\nabla d(x,x_0), \nabla V\> = \<\nabla V, U\>(x_0)+
\int_0^{d(x,x_0)} Hess_V(U,U) (x_s) ds \ge  \delta d(x,x_0)-c_1.
$$
So condition (\ref{cor2b}) holds. Furthermore if $Hess_V\ge
\delta>(1+\sqrt{2})\sigma \sqrt{d-1}$, Wang \cite{Wang08} proves the
LSI for $\mu$. When $\sigma \sqrt{d-1}<\delta\le (1+\sqrt{2})\sigma
\sqrt{d-1}$, the LSI is actually unknown. Also see \cite{ATW} for
the Harnack type inequality on this type of manifold.

One main feature of our condition (\ref{cor2b}) is: it demands only
on the radial derivative of $V$, NOT on $Hess_V$.

\end{remark}

\begin{proof} At first we borrow the proof of \cite[Lemma
2.1]{Wang08} for controlling $\Delta \rho$ where $\rho(x)=d(x,x_0)$.
By (\ref{cor2a}) and the Laplacian comparison theorem, we have for
$x\notin cut(x_0)$ different from $x_0$
$$
\Delta \rho \le \sqrt{(c+\sigma^2 \rho^2)(d-1)} \coth \left(
\rho\sqrt{(c+\sigma^2 \rho^2)/(d-1)} \right).
$$
Then outside of $cut(x_0)$ we get
$$\aligned
\Delta \rho^2&=2\rho \Delta \rho +2\\
&\le 2\rho  \sqrt{(c+\sigma^2 \rho^2)(d-1)} \coth \left(
\rho\sqrt{(c+\sigma^2 \rho^2)/(d-1)} \right)+2\\
&\le 2d + 2\rho  \sqrt{(c+\sigma^2 \rho^2)(d-1)}\endaligned
$$
where the last inequality follows by $r\cosh r\le (1+r)\sinh r\
(r\ge0)$. It is well known that $\Delta \rho$ in the distribution
sense gives a non-positive measure on $cut(x_0)$, the above
inequality holds in the distribution sense over $E$.

Hence under the condition that $\delta>\sigma \sqrt{d-1}$, for
$U=\lambda \rho^2$ where $0<\lambda< \frac 12(\delta - \sigma
\sqrt{d-1})$, we have in the sense of distribution
$$\aligned
\LL U + |\nabla U|^2 &\le 2\lambda [2d + 2\rho  \sqrt{(c+\sigma^2
\rho^2)(d-1)}]- 2\lambda \rho \<\nabla \rho, \nabla V\> + 4\lambda^2
\rho^2\\
&\le - c \rho^2 +b
\endaligned$$
for some positive constants $b,c$, i.e. condition (\ref{lyapT2b}) is
satisfied. So the $W_2H$ inequality follows by Theorem
\ref{thm1}(1).
\end{proof}

Our condition ``$\delta>\sigma \sqrt{d-1}$" for $W_2H$ is sharp as
shown by the following example taken from \cite{Wang08}.

\begin{example}\label{exam2} {\rm Let $E=\R^2$ be equipped with the following
Riemannian metric
$$
ds^2=dr^2 + (re^{kr^2}) d\theta^2$$ under the polar coordinates
$(r,\theta)$, where $k>0$ is constant. Then $Ric_{(r,\theta)}=-4k -
4k^2 r^2$. Then (\ref{cor2a}) holds with $\sigma=2k$. Let
$V:=\frac\delta 2 r^2$, which satisfies (\ref{cor2b}). If
$\delta>\sigma\sqrt{d-1}=2k$, we have $W_2H$. But if $\delta\le
\sigma\sqrt{d-1}=2k$, $e^{-V} dx=r e^{kr^2-\delta r^2/2} dr d\theta$
is infinite measure, so that $W_2H$ does not hold.

 }
\end{example}

\section{Proof of Theorem \ref{thm1}}

\subsection{Several lemmas}

As was recalled in a previous remark, we may assume without loss of
generality that $\mu$ verifies a Poincar\'e inequality with constant
$C_P$, i.e. $\int g^2 d\mu\le C_P \int |\nabla g|^2d\mu$ for all
smooth $g$ with $\mu(g)=0$.

We begin with the following
\begin{lemma} \label{lem21}{\rm (\cite[Theorem 1.13]{CG05})} If $\mu$ satisfies the Poincar\'e inequality with
constant $C_P$, then for all smooth and bounded $g$,
\begin{equation}
Ent_\mu(g^2)\le 2C_P \left(2\log2+\frac12\log
\frac{\|g^2\|_\infty}{\mu(g^2)}\right) \int |\nabla g|^2d\mu.
\end{equation}
Conversely, if the preceding restricted logarithmic Sobolev is true then $\mu$ satisfies a Poincar\'e inequality with constant $4C_P\log2$.
\end{lemma}

\begin{lemma}\label{lem24} Assume that the following restricted logarithmic Sobolev inequality holds: there
exist constants $\eta, C_{\eta}>0$ such that $$\Ent_\mu(g^2)\le
2C_{\eta} \int |\nabla g|^2d\mu$$ for all smooth and bounded
functions $g$ satisfying
\begin{equation}\label{restriction}
g^2\le\left(\int g^2d\mu\right) e^{2\eta (d^2(x,x_0)+\int
d^2(y,x_0)d\mu(y))}.\end{equation} Then $\mu$ satisfies $W_2H(C)$
with  $C=\max\{C_\eta; (2\eta)^{-1}\}$.
\end{lemma}
\begin{proof}  We recall the (short
and simple) proof from \cite[Theorem 1.17]{CG05}.

Given a fixed bounded $f$ with $\mu(f)=0$ consider for any
$\lambda\in\R$, $g^2_\lambda:=e^{\tilde \eta Q(\lambda f)}$ where
$\tilde \eta:=\min\{1/(2C_\eta);\eta\}\in (0,\eta]$. By the
definition of $Q$ we easily get
$$Q(\lambda f)(x)\le \int (\lambda f(y) + d^2(x,y)) d\mu(y)\le 2d^2(x,x_0)+2\int d^2(y,x_0)\mu(dy).$$
Let $G(\lambda)=\mu(g^2_\lambda)$. By  Bobkov-Goetze's criterion (Theorem 1.4(2)), if $G(1)\le 1$
(for all such $f$), then $W_2H(C)$ holds with $C=1/(2\tilde \eta)=\max\{C_\eta; (2\eta)^{-1}\}$.
Assume by absurd that $G(1)>1$. Introduce $\lambda_0=\inf \{\lambda\in[0,1];~G(u)>1,\forall
u\ge\lambda\}$, and remark that $\lambda_0<1$, $G(\lambda_0)=1$ as well as $G(0)=1$ and that
$G(\lambda)>1$ as soon as $\lambda\in]\lambda_0,1]$.

Note at first that if $G(\lambda)\ge 1$ then
$$g_\lambda^2\le e^{2\tilde \eta(d^2(x,x_0) + \int d^2(x,x_0)d\mu(x))}\le G(\lambda)e^{2\eta
(d^2(x,x_0) + \int d^2(x,x_0)d\mu(x))}$$ i.e., $g_\lambda$ satisfies  condition
(\ref{restriction}). Since $Q_tf(x):=\inf_{y\in E}(f(y)+ \frac 1{2t} d^2 (x,y))$ is the Hopf-Lax
solution of the Hamilton-Jacobi equation: $\partial_t Q_t f + \frac 12 |\nabla Q_tf|^2=0$
(\cite{BGL01}) and $Q(\lambda f)=\lambda Q_{\lambda/2} f$, we have
$$\lambda G'(\lambda)=\int g^2_\lambda\log g^2_\lambda d\mu-\frac1{\tilde \eta}\int|\nabla g_\lambda|^2d\mu.$$

Since $\tilde \eta=\min\{1/(2C_\eta);\eta\}$, the restricted
logarithmic Sobolev inequality in Lemma \ref{lem24} yields for
$\lambda\in]\lambda_0,1]$
$$\lambda G'(\lambda)\le G(\lambda)\log G(\lambda)$$
which is nothing else than the differential inequality
$(\lambda^{-1}\log G(\lambda))'\le0$. That implies that
$\lambda^{-1}\log G(\lambda)$ is nonincreasing so that
$$\log G(1)\le \frac{\log(G(\lambda_0))}{\lambda_0}$$
(taken as limit $\lim_{\lambda\to
0}\frac{\log(G(\lambda))}{\lambda}=0$ if $\lambda_0=0$). It readily
implies that $G(1)\le1$ which is the Bobkov-Goetze's condition.
\end{proof}

\begin{remark}
The fact that the restricted logarithmic Sobolev inequality implies
$W_2H$ inequality was proven in \cite[Th. 1.17]{CG05}. In addition a
Hardy criterion
for this inequality on the real line is given in \cite[Prop. 5.5]{CG05}.
\end{remark}

Let $(\EE, \D(\EE))$ be the Dirichlet form associated with $\LL$ in
$L^2(\mu)$. It is the closure of $\EE(f,g)=\<-\LL f,
g\>_{L^2(\mu)}=\int \nabla f\cdot \nabla g d\mu, \ f,g\in
C_0^\infty(E)$ by the essential self-adjointness of $(\LL,
C_0^\infty(E))$.

\begin{lemma}\label{lem23} Let $U$ be a nonnegative locally Lipschitzian function such that $\LL U +|\nabla U|^2\le
-\phi$ in the distribution sense, where $\phi$ is lower bounded,
then for any $g\in \D(\EE)$,
\begin{equation}\label{lem23e}
\int \phi g^2 d\mu \le \EE(g,g).
\end{equation}
\end{lemma}

\begin{proof} As $\phi\wedge N$ satisfies also the condition,
if (\ref{lem23e}) is true with $\phi\wedge N$, then it is true with
$\phi$ by letting $N\to+\infty$. In other words we can and will
assume that $\phi$ is bounded.

One can approach any $g\in \D(\EE)$ by $(g_n)\subset C_0^\infty(E)$
: $\int (g_n-g)^2 d\mu + \EE(g_n-g, g_n-g)\to 0$. Thus is enough to
prove (\ref{lem23e}) for $g\in C_0^\infty(E)$. For $g\in
C_0^\infty(E)$, we have
$$\int (-\LL U) g^2 d\mu= \int U(-\LL g^2) d\mu=\int \nabla U\cdot
\nabla (g^2) d\mu$$ where the first equality comes from the
definition of the distribution $-\LL U$ and a direct calculus, the
second one is true at first for $U\in C_0^\infty(E)$ and is extended
at first to any Lipschitzian $U$ with compact support, then to any
locally Lipschitzian $U$.

Thus using $2  g \nabla U \cdot \nabla g\le |\nabla U|^2g^2+|\nabla
g|^2$, we get
$$\aligned
\int \phi g^2 d\mu &\le \int (-\LL U - |\nabla U|^2) g^2
d\mu\\
&=\int  \left(2 g \nabla U \cdot \nabla g - |\nabla U|^2 g^2
\right)d\mu \le \int |\nabla g|^2 d\mu
\endaligned$$
which is the desired result.
\end{proof}

We also require the consequence below of the Lyapunov condition
(\ref{lyapT2b}).

\begin{lemma}\label{lem22} If the Lyapunov condition (\ref{lyapT2b}) holds, then
there exist $\delta>0$, $x_0\in E$ such that
\begin{equation}\label{expdcarre}
\int e^{\delta d^2(x,x_0)}d\mu<\infty.
\end{equation}
\end{lemma}
\begin{proof} Under the condition (\ref{lyapT2b}), $\LL$ satisfies a spectral
gap property in $L^2(\mu)$ and then by \cite{GLWY07}, the following $W_1I$-inequality holds:
$$
W_1^2(\nu,\mu) \le 4C^2 I(\nu|\mu), \ \forall \nu\in M_1(E)
$$
where
\begin{equation}\label{I}I(\nu|\mu):=\begin{cases}\EE(\sqrt{h},\sqrt{h}),
\ &\text{ if }\  \nu=h\mu, \ \sqrt{h}\in\D(\EE)\\
+\infty, &\text{ otherwise }
\end{cases}\end{equation} is the so called Fisher
information. By \cite{GLWW08}, the above $W_1I$-inequality is
stronger than $W_1H(C)$, which is equivalent to the gaussian
integrability (\ref{expdcarre}).
\end{proof}

It would be interesting to find a simple or direct argument leading
to (\ref{expdcarre}).

\subsection{Proof of Theorem \ref{thm1}(1)}
Choose $\eta>0$ such that $\eta<\min(1,\delta/2)$ where $\delta$
comes from the gaussian integrability condition (\ref{expdcarre})
which holds by Lemma \ref{lem22}. We have only to prove the
restricted LSI in Lemma \ref{lem24} under the Lyapunov condition
(\ref{lyapT2}).

To simplify the notation, define $M=e^{2\eta \int d^2(x,x_0)d\mu(x)}$. \\
Let $h=g^2$ be positive and smooth  with $\mu(h)=1$ and $h\le M
e^{2\eta d^2(x,x_0)}$. By Lemma \ref{lem22} and our choice of
$\eta$, $\int h\log hd\mu$ is bounded by some constant, say
$c(\eta,\mu)$. Take $K>e$, to be chosen later. We have
\begin{equation}\label{lem23b}\aligned
\int h\log hd\mu&= \int_{h\le K} h\log hd\mu+\int_{h> K} h\log hd\mu\\
&\le\int (h\wedge K)\log (h\wedge K)d\mu+(\log M) \int_{h> K} hd\mu+
2\eta \int_{h> K} hd^2(x,x_0)d\mu.\endaligned
\end{equation}
As $\int_{h\le K} h\log hd\mu\ge \int_{h\le K} (h-1)d\mu\ge -
\int_{h>K} hd\mu $, we have
$$\int h\log hd\mu\ge \int_{h>K}h\log hd\mu-\int_{h>K}hd\mu.$$
It yields
$$\int_{h>K} hd\mu\le \frac1{\log K}\int_{h>K} h\log hd\mu\le \frac{1}{\log K}\left(\int h\log hd\mu +
\int_{h>K}hd\mu\right)$$ so that
\begin{equation}\label{lem23c}\int_{h>K} hd\mu\le \frac1{\log K-1}\int h\log
hd\mu\le \frac{c(\eta,\mu)}{\log K-1}.\end{equation}

\eqref{lem23c} furnishes an immediate useful bound for the second term in the right hand side of
(\ref{lem23b}). Indeed, if $3\log M\le \log K -1$ then
$$ \log M\int_{h> K} h d\mu\le \frac13 \int h\log hd\mu.$$
Remark also that for $K>e$
$$1\ge \int h\wedge K d\mu\ge 1-\frac{c(\eta,\mu)}{\log K-1}$$
so that for $K$ large enough  (independent of $h$), $\int h\wedge Kd\mu\ge 1/2$ and thus by Lemma
\ref{lem21}
$$\aligned
\int(h\wedge K)\log (h\wedge K)d\mu &\le \int(h\wedge K)\log
\left(\frac{h\wedge K} {\int h\wedge Kd\mu}\right)d\mu\\
& \le C_P(2\log2+\frac12\log (2K))\int |\nabla \sqrt{h}|^2 d\mu.\endaligned$$ We then only have to
bound the last term in \eqref{lem23b}. Unfortunately, we cannot  directly apply the Lyapunov
condition due to a lack of regularity of $h\BBone_{h>K}$. So we first regularize this function. To
this end, introduce  the map $\psi$ with
$$\psi(u)=\left\{\begin{array}{ll}
0&{\rm if~}0\le u\le \sqrt{K/2}\\
\frac{\sqrt{2}}{\sqrt{2}-1}(u-\sqrt{K/2})&{\rm if~}\sqrt{K/2}\le u\le\sqrt{K}\\
u&{\rm if~}\sqrt{K}\le u. \end{array}\right.$$  Now using Lyapunov
condition (\ref{lyapT2b}) and Lemma \ref{lem23} (applicable for
$\psi(\sqrt{h})=\psi(g)$ is locally Lipschitzian), we have
\begin{eqnarray*}
2\eta\int_{h>K}hd^2(x,x_0)d\mu&\le& 2\eta \int \psi^2(\sqrt{h})d^2(x,x_0)d\mu\\
&\le& \frac{2\eta}{c} \int \psi^2(\sqrt{h})[c d^2(x,x_0) -b]d\mu
+\frac{2\eta b}{c}\int \psi^2(\sqrt{h}) d\mu\\
&\le & \frac{2\eta}{c} \int  |\nabla \psi(\sqrt{h})|^2 d\mu
+\frac{2\eta b}{c}\int \psi^2(\sqrt{h}) d\mu\\
&\le & \frac{4\eta}{c(\sqrt{2} -1)^2}\int  |\nabla \sqrt{h}|^2 d\mu
+\frac{2\eta b}{c}\int \psi^2(\sqrt{h}) d\mu .
\end{eqnarray*}
As $\psi^2(\sqrt{h})\le h1_{h>K/2}$, the lat term above can be
bounded by $(1/3)\int h\log h d\mu$ if $K$ is large enough so that
$2\eta b c^{-1}\le (\log(K/2)-1)/3$, by (\ref{lem23c}).

Plugging all those estimates into (\ref{lem23b}), we obtain the
desired restricted LSI.

\subsection{Proof of Theorem \ref{thm1}(2)}
Our argument will be a combination of the Lyapunov condition,
leading to defective $W_2I$ inequality and the HWI inequality of
Otto-Villani.

We begin with the following fact (\cite[Proposition 7.10]{Vill03}):
\begin{equation}\label{eqenplus}
W_2^2(\nu,\mu) \le 2 \|d(\cdot, x_0)^2(\nu-\mu)\|_{TV}.
\end{equation}
Now for every function $g$ with $|g|\le \phi(x):=c d(x,x_0)^2$, we
have by (\ref{lyapT2b}) and Lemma \ref{lem23},

$$
\aligned
\int g d(\nu-\mu) &\le \nu(\phi) + \mu(\phi)\\
&\le \int \left(- c d^2(x,x_0) + b\right) d\nu(x) + \mu(\phi)\\
&\le I(\nu|\mu) + b + \mu(\phi)
\endaligned
$$
 Taking the
supremum over all such $g$, we get

$$
\frac {c}2 W_2^2(\nu,\mu) \le c \|d(\cdot, x_0)^2(\nu-\mu)\|_{TV} \le I(\nu|\mu) + b + \mu(\phi),
$$
which yields thanks to \eqref{eqenplus}
$$
W_2^2(\nu,\mu) \le \frac 2c I(\nu|\mu) + \frac 2c[b + \mu(\phi)].
$$
Substituting it into the HWI inequality of Otto-Villani
\cite{OVill00} (or for  its Riemannian version by
Bobkov-Gentil-Ledoux \cite{BGL01}):
\begin{equation}\label{HWI} H(\nu|\mu)\le 2\sqrt{I(\nu|\mu)}W_2(\nu,\mu)-\frac K2W_2^2(\nu,\mu),\end{equation}
and using $2ab \le \varepsilon a^2 + \frac 1 \varepsilon b^2$ we finally get
\begin{equation}\label{thm52b}
\aligned
H(\nu|\mu)
&\le \varepsilon I(\nu|\mu) +\left(1 - \frac K2 + \frac 1\varepsilon\right)W_2(\nu,\mu)^2\\
&\le A I(\nu|\mu) +B
\endaligned
\end{equation}
where
$$A=(1-\frac K2)\frac 2c + \varepsilon,\ \ B=\frac 2c[b + \mu(\phi)]\left(1 - \frac K2 + \frac 1\varepsilon\right). $$
This inequality is sometimes called a defective log-Sobolev inequality. But it is well known by
Rothaus' lemma, that a defective log-Sobolev inequality together with the spectral gap implies the
(tight) log-Sobolev inequality
$$
H(\nu|\mu) \le [A + (B+2)C_P)] I(\nu|\mu).
$$
The proof is completed.

\begin{remark} {\rm If for any $c>0$, there are $U,b$ such that the Lyap condition (\ref{lyapT2b}) holds, then the defective LSI (\ref{thm52b}) becomes the so called {\it super-LSI}, which is equivalent to the supercontractivity
of
the semigroup $(P_t)$ generated by $\LL$, i.e., $\|P_t\|_{L^p\to
L^q}<+\infty$ for any $t>0, q>p>1$.
}
\end{remark}

\section{Some further remarks}
\subsection{A generalization of Corollary \ref{cor2}}

\begin{corollary} Assume that
$$
Ric_x \ge -\alpha(d(x,x_0))
$$
where $\alpha(r)$ is some positive increasing function on $\R^+$,
and
$$
\<\nabla d(x,x_0), \nabla V\> \ge \beta(d(x,x_0)) -b
$$
for some constant $b>0$ and some positive increasing function
$\beta$ on $\R^+$. If
\begin{equation}\label{cor3c} \beta(r) - \sqrt{\alpha(r)}\ge \eta r, \ r>0  \end{equation}
for some positive constant $\eta$,
 then $\mu$ satisfies $W_2H$.
\end{corollary}

\begin{proof}As in the proof recalled in Corollary \ref{cor2}, for $\rho=d(x,x_0)$,
by the Laplacian comparison theorem, there is some constant $c_1>0$
such that
$$
\Delta \rho^2 \le c_1(1+\rho) + 2\rho\sqrt{\alpha(\rho)}
$$
at first outside of $cut(x_0)$ then in distribution over $E$.
Consequently by condition (\ref{cor3c}) there are positive constants
$c_2<2\eta, c_3$ such that
$$
\LL \rho^2 =\Delta \rho^2 - 2\rho\<\nabla \rho, \nabla V\>\le
c_1(1+\rho) +2\rho(\sqrt{\alpha(\rho)}- \beta(\rho)+b)\le -c_2\rho^2
+c_3.
$$
Now for $U=\lambda \rho^2$, it satisfies (\ref{lyapT2b}) when
$\lambda>0$ is small enough. Then the $W_2H$ follows by Theorem
\ref{thm1}.
\end{proof}

\subsection{LSI in the unbounded curvature case}
We now generalize the LSI in Theorem \ref{thm1} in the case where
Bakry-Emery's curvature  is not lower bounded, by means of the
above-tangent lemma.

\begin{proposition}\label{prop41} Assume that
\begin{equation}\label{prop41a}
Ric_x + Hess_V \ge -\Phi(d(x,x_0)) \end{equation} where $\Phi$ is
some positive non-decreasing continuous function on $\R^+$, and
there is some nonnegative locally Lipschitzian function $U$ such
that for some constants $b,c>0$
\begin{equation}\label{prop41b}
\LL U + |\nabla U|^2 \le - c d^2(x,x_0)\Phi(2d(x,x_0)) +b
\end{equation}
in distribution, then $\mu$ satisfies the LSI.
\end{proposition}

\begin{proof} Instead of the HWI in the proof of the LSI in Theorem \ref{thm1}, we
go back to the above-tangent lemma (see \cite[Theorem 7.1]{BK07} and
references therein) : for two probability measures $\nu=h\mu, \tilde
\nu=\tilde h\mu$ with smooth and compactly supported densities $h,
\tilde h$, let $T(x):=\exp_x(\nabla \theta)$ (where $\theta$ is some
``convex" function) be the optimal transport pushing forward $\nu$
to $\tilde\nu$ and realizing $W_2^2(\nu,\tilde \nu)$. Then
\begin{equation}\label{prop41c}
Ent_\mu(h) \le Ent_\mu(\tilde h)- \int \<\nabla \theta, \nabla h\>
d\mu + \int \DD_V(x, T(x)) h d\mu
\end{equation}
where $\DD_V(x,T(x))$ is the defect of the convexity of $V$, defined
by
$$
\DD_V(x,T(x))=-\int_0^1 (1-t) \left(Ric_{\gamma(t)}+ Hess_{V,
\gamma(t)} \right)(\dot\gamma(t), \dot\gamma(t)) dt.
$$
Here $\gamma(t)=\exp_x(t\nabla \theta)$ is the geodesic joining $x$
to $T(x)$.

Choose a sequence of $\mu$-probability measures $\mu_n:=h_n\mu$ with
$h_n \in C_0^\infty(E)$, such that $W_2(\mu_n,\mu)\to 0$ and
$I(\mu_n|\mu)\to 0$ (recalling that the condition (\ref{prop41b}),
stronger than (\ref{lyapT2b}), implies the Gaussian integrability of
$\mu$ by Lemma \ref{lem22}). Below we apply the above-tangent lemma
to $(\nu, \tilde \nu=\mu_n)$

The first term on the right hand of (\ref{prop41c}) is easy to
control by Cauchy-Schwarz:
$$
|\int \<\nabla \theta, \nabla h\> d\mu|= |\int 2 \sqrt{h}\<\nabla
\theta, \nabla \sqrt{h}\> d\mu|\le 2\sqrt{\int |\nabla \theta|^2 h
d\mu \int |\nabla \sqrt{h}|^2 d\mu} = 2 W_2(\nu,\mu_n)\sqrt{
I(\nu|\mu)}.
$$
Now we treat the last term in (\ref{prop41c}). By our condition,
$$
\DD_V(x,T(x)) \le \int_0^1(1-t) \Phi(d(\gamma(t), x_0)) |\nabla
\theta|^2 dt.
$$
Note that $|\nabla \theta|=d(x,T(x))\le 2\max\{d(x, x_0),
d(T(x),x_0)\}$ and using $d(\gamma(t), x_0)\le d(x,x_0) +
td(x,T(x))$ for $t\in [0,1/2]$ and $d(\gamma(t), x_0)\le d(T(x),x_0)
+ (1-t)d(x,T(x))$ for $t\in [1/2,1]$, $d(\gamma(t), x_0)\le
2\max\{d(x, x_0), d(T(x),x_0)\}$. We thus obtain
$$\aligned
\int \DD_V(x, T(x)) h d\mu&\le 2\int \Phi(2\max\{d(x, x_0),
d(T(x),x_0)\})\max\{d(x, x_0)^2, d(T(x),x_0)^2\} hd\mu\\
&\le 2\left(\int \Phi(2d(x, x_0))d(x, x_0)^2 h d\mu + \int
\Phi(2d(T(x), x_0))d(T(x), x_0)^2 h d\mu \right)\endaligned
$$
By Lemma \ref{lem23} and our condition (\ref{prop41b}),
$$\aligned
&c \int \Phi(2d(x, x_0))d(x, x_0)^2 h d\mu\le b + I(\nu|\mu)\\
&c \int \Phi(2d(T(x), x_0))d(T(x), x_0)^2 h d\mu \le b +
I(\mu_n|\mu)\endaligned$$ Plugging those estimates into
(\ref{prop41c}) and letting $n\to\infty$, we get finally
$$
H(\nu|\mu) \le 2W_2(\nu, \mu)\sqrt{I(\nu|\mu)} + \frac1c(I(\nu|\mu)
+2b)
$$
Again using Lemma \ref{lem23} and our condition (\ref{prop41b}), we
have
$$
W_2^2(\nu,\mu)\le 2 \left(\int d^2(x,x_0) d\mu + \int d^2(x,x_0)
d\mu\right)\le \frac2{c\Phi(0)}(I(\nu|\mu) +2b).
$$
Consequently we obtain the defective LSI:
\begin{equation}\label{prop41d}
H(\nu|\mu) \le 2\sqrt{\frac 2{c\Phi(0)}} (I(\nu|\mu)
+b)+\frac1c(I(\nu|\mu) +2b)
\end{equation}
where the LSI follows for the spectral gap exists under
(\ref{prop41b}).
\end{proof}

\begin{remark} Under (\ref{prop41a}), if for any $c>0$ there are $U,b$ such that the Lyapunov function
condition (\ref{prop41b}) holds, the defective LSI
(\ref{prop41d}) says that for any $\varepsilon>0$, there is some
constant $B(\varepsilon)$ such that
$$
H(\nu|\mu)\le \varepsilon I(\nu|\mu) + B(\varepsilon), \ \nu\in
M_1(E)
$$
which is well known to be equivalent to the supercontractivity of
the semigroup $(P_t)$ generated by $\LL$, i.e., $\|P_t\|_{L^p\to
L^q}<+\infty$ for any $t>0, q>p>1$.
\end{remark}

\begin{remark} Barthe and Kolesnikov \cite{BK07} used the
above-tangent lemma to derive  modified LSI and isoperimetric
inequalities. One aspect of their method consists in controlling the
defective term $\int \DD_V(x,T(x)) h d\mu$ by $c Ent_\mu(h) +b$ for
some positive constant $c<1$, by using some integrability condition
on $\mu$ (as in Wang's criterion). Our method is to bound that
defective term by  $c I(\nu|\mu) +b$, by means of the Lyapunov
function: the advantage here is that constant $c>0$ can be
arbitrary.
\end{remark}

\begin{example}{\rm Let $E=\R^2$ equipped with the Euclidean metric.
For any $p> 2$ fixed, consider $V=r^p(2+\sin (k\theta))$, where
$(r,\theta)$ is the polar coordinates system and $k\in\N^*$. Since
$$
\Delta V = r^{p-2}[p^2 (2+\sin (k\theta))-k^2 \sin (k\theta) ]
$$
Assume $k> \sqrt{3}p$. Then in the direction $\theta$ such that
$\sin(k\theta)=1$, $Hess_V \le - \frac 12 (k^2-3p^2) r^{p-2}$, i.e.,
the Bakry-Emery curvature is very negative  and no known result
exists in such case.

It is easy to see that condition (\ref{prop41a}) is verified with
$\Phi(r)=a r^{p-2}$ for some $a>0$. Taking $U= r^2$, we see that
$$
\LL U + |\nabla U|^2 = 4  - 2 p r^p (2+\sin(k\theta))+4r^2
$$
i.e., condition (\ref{prop41b}) is satisfied. We get thus the LSI
for $\mu$ by Proposition \ref{prop41}. }\end{example}

\subsection{A Lyapunov condition for Gozlan's weighted Poincar\'e inequality}

 As mentionned before, in a recent work, Gozlan \cite{Goz07b} proved that $W_2H$ inequality on $E=\R^d$ is
 implied by a weighted Poincar\'e inequality
 $$\Var_\mu(f)\le c \int \sum_1^d\frac1{1+x_i^2}\left(\frac{\partial f}{\partial x_i}\right)^2d\mu.$$
  In dimension one, a
 Hardy criterion is available for this weighted Poincar\'e inequality
 which is not the same as the one from \cite{CG05}.
 Note however that
 this weighted Poincar\'e inequality, as stronger than Poincar\'e inequality,
 can be shown to imply a converse weighted Poincar\'e
 inequality (the weight is now in the variance), by a simple change of function argument, and in dimension one a Hardy's
 criterion is also available for this inequality which is in fact the
 same as the one for the restricted logarithmic Sobolev inequality.\\
From this, we conclude that in fact, in the real line case, the
restricted logarithmic Sobolev inequality is in fact implied by Gozlan's weighted Poincar\'e inequality.
Whether it is the case in any dimension would have to be investigated.

It is however quite easy, following \cite{BBCG} to give a Lyapunov condition for Gozlan's weighted Poincar\'e inequality on $\R^d$.

\begin{theorem}\label{thm4-11}
Let $w_i=w_i(x_1,...,x^d)$ be positive for all $(x_1,...,x_d)\in
\R^d$, and $\omega_i>\epsilon_r>0$ on $B(0,r)$. Introduce the
diffusion generator
$$\widetilde\LL  = \sum_{i=1}^d\left( \omega_i\partial^2_{i}+(\partial_i\omega_i -\omega_i\partial_i V)\partial_i\right),$$
where $\partial_i=\partial/\partial x_i$. Suppose now that there exists $W\ge1$, $\lambda,b>0$ and $R>$ such that
\begin{equation}\label{lyapSGW}
\widetilde\LL W\le -\lambda W+b\BBone_{B(0,R)}
\end{equation}
then $\mu$ verifies a weighted Poincar\'e inequality with some
constant $c>0$
\begin{equation}
\label{SGW}
\Var_\mu(f)\le c\int\sum_{i=1}^d\omega_i\left(\partial_i f\right)^2d\mu.
\end{equation}
\end{theorem}
\begin{proof}
The proof follows exactly the line of the one of \cite{BBCG} once it has been remarked that
\begin{enumerate}
\item $\widetilde L$ is associated to the Dirichlet form $\widetilde\EE(f,g)=-\int f\widetilde\LL gd\mu$ on $L^2(\mu)$,
 reversible w.r.t. $\mu$ and $\widetilde \EE(f,f)=\int\sum_{i=1}^d\omega_i\left(\partial_i f\right)^2d\mu$;
\item a local weighted Poincar\'e inequality is valid for this Dirichlet form as $\omega_i>\epsilon_r>0$ on $B(0,r)$ (as a local Poincar\'e inequality is available on balls).
\end{enumerate}
\end{proof}
\begin{remark}
Our setting is a little bit more general than Gozlan \cite{Goz07b} concerning the assumption on $\omega$ but with the additional term $\partial_i\omega_i\partial_i W$ in the sum. Note once again that they are a little bit more difficult to handle than the one in Corollary \ref{cor1} and still not comparable. \\
One of the major points of Gozlan's weighted Poincar\'e inequality is, in the case where $\omega_i(x_1,...,x_d)=\omega(x_i)$, in fact equivalent to some transportation-information inequality (with an unusual distance function) when $\omega$ satisfies some conditions (namely, $\omega=\sqrt{\tilde\omega'}$ where $\tilde\omega$ is odd, at least linearly increasing). However, when $\omega_i=1/(1+x_i^2)$, this transportation inequality is stronger than $W_2H$.
\end{remark}
We end up this note with some final conditions ensuring $W_2H$,
similar to Gozlan's one (see Remark 2.4(2)).
\begin{corollary}
In the setting of Corollary \ref{cor1}.\item 1) If there are
positive constants  $a<1, R,c>0,$ such that
\begin{equation}
\label{kustrWSP} \sum_{i=1}^d\left((1-a)\omega_i(\partial_i
V)^2-\partial_i\omega_i\partial_iV-\omega_i\partial^2_{i}V\right)\ge
c,  \forall x: |x|>R
\end{equation}
is verified then the weighted Poincar\'e inequality (\ref{SGW}) is verified.\\
2) In particular, consider $\omega_i(x_1,...,x_d)=(1+x_i^2)^{-1}$,
if there are positive constants  $a<1, R,c>0,$ such that for all
$x\in \R^d$ with $|x|>R$, one of
\begin{equation}
\label{kustrWSP2}\sum_{i=1}^d\left((1-a)(\partial_i
V)^2+\frac{2x_i}{1+x_i^2}\partial_iV-\partial^2_{i}V\right)\frac1{1+x_i^2}\ge
c
\end{equation}
or
\begin{equation}
\label{simplWSG} \sum_{i=1}^d \left(\frac{x_i\partial_i
V}{1+x_i^2}-\frac{1-x_i^2}{(1+x_i^2)^2}\right)\ge c
\end{equation}
is verified, then $W_2H$ holds.
\end{corollary}
\begin{proof} Part 1) is a particular case of Theorem \ref{thm4-11}
with $W=e^{aV}$, together with Gozlan's result. Condition
(\ref{kustrWSP2}) is just a particular version of part 1). The last
case under condition (\ref{simplWSG}) comes from Theorem
\ref{thm4-11} with $W=e^{a|x|^2}$ for sufficiently small $a$.
\end{proof}


\bibliographystyle{plain}

\end{document}